\documentstyle{article}
\begin{document}

\newtheorem{thm}{Theorem}
\newtheorem{defin}[thm]{Definition}
\newtheorem{lemma}[thm]{Lemma}
\newtheorem{propo}[thm]{Proposition}
\newtheorem{cor}[thm]{Corollary}
\newtheorem{conj}[thm]{Conjecture}

\centerline{\huge \bf The Area of the Medial Parallelogram}

\centerline{\huge \bf of a Tetrahedron}
\bigskip
\smallskip

\centerline{\parbox{56mm}{David N. Yetter \\ Department of
Mathematics \\ Kansas State University \\ Manhattan, KS 66506}
\footnote{Supported by NSF Grant \# DMS-9504423}}
\bigskip
\bigskip

	The purpose of this note is to present a simple formula for the
area of the medial parallelogram of a tetrahedron in terms of the lengths
of the six edges, together with its derivation.  It is not claimed that
this formula is original, though diligent searches of standard sources of 
mensuration formulae, and inquiries to mathematicians in a variety of
fields suggest that it is not generally known, and {\em may} be new. The author
would be grateful for any references to previously published instances
of this formula.

	Despite the very classical nature of the problem it solves, there
is some serious contemporary interest arising from recently proposed
simplicial models for quantum gravity (cf. \cite{Barbieri}, 
\cite{Barrett-Crane}), in which such a formula is needed to approach the
problem of length operators.

\begin{figure}

\setlength{\unitlength}{0.00083300in}%
\begingroup\makeatletter\ifx\SetFigFont\undefined
% extract first six characters in \fmtname
\def\x#1#2#3#4#5#6#7\relax{\def\x{#1#2#3#4#5#6}}%
\expandafter\x\fmtname xxxxxx\relax \def\y{splain}%
\ifx\x\y   % LaTeX or SliTeX?
\gdef\SetFigFont#1#2#3{%
  \ifnum #1<17\tiny\else \ifnum #1<20\small\else
  \ifnum #1<24\normalsize\else \ifnum #1<29\large\else
  \ifnum #1<34\Large\else \ifnum #1<41\LARGE\else
     \huge\fi\fi\fi\fi\fi\fi
  \csname #3\endcsname}%
\else
\gdef\SetFigFont#1#2#3{\begingroup
  \count@#1\relax \ifnum 25<\count@\count@25\fi
  \def\x{\endgroup\@setsize\SetFigFont{#2pt}}%
  \expandafter\x
    \csname \romannumeral\the\count@ pt\expandafter\endcsname
    \csname @\romannumeral\the\count@ pt\endcsname
  \csname #3\endcsname}%
\fi
\fi\endgroup
\begin{picture}(4212,4389)(3139,-6088)
\thicklines
\put(3151,-5761){\line( 0, 1){  0}}
\put(3151,-5761){\line( 1, 0){4050}}
\put(7201,-5761){\line(-1, 3){1350}}
\put(5851,-1711){\line(-2,-3){2700}}
\put(5851,-1711){\line( 2,-3){1350}}
\put(7201,-3736){\line(-2,-1){4050}}
\put(7201,-3736){\line( 0,-1){2025}}
\put(4276,-3511){\makebox(0,0)[lb]{\smash{\SetFigFont{12}{14.4}{it}a}}}
\put(6301,-3736){\makebox(0,0)[lb]{\smash{\SetFigFont{12}{14.4}{it}b}}}
\put(7351,-4861){\makebox(0,0)[lb]{\smash{\SetFigFont{12}{14.4}{it}c}}}
\put(5101,-6061){\makebox(0,0)[lb]{\smash{\SetFigFont{12}{14.4}{it}d}}}
\put(6676,-2611){\makebox(0,0)[lb]{\smash{\SetFigFont{12}{14.4}{it}e}}}
\put(5026,-4636){\makebox(0,0)[lb]{\smash{\SetFigFont{12}{14.4}{it}f}}}
\end{picture}

\caption{\label{tet} A generic tetrahedron}
\end{figure}
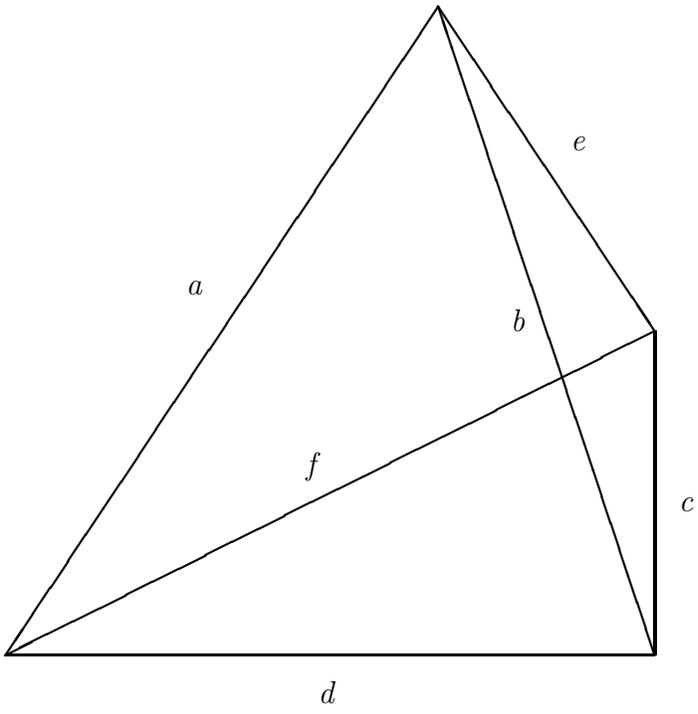

	Consider a tetrahedron with edge-lengths as in Figure \ref{tet}.
Fix a pair of non-incident edges, say those of lengths $d$ and $e$.  
It is then easy to see that the midpoints of the remaining four edges
lie in a plane parallel to both of the chosen edges, and equidistant from
the planes containing each chosen edge and parallel to both, and that they
form the vertices of a parallelogram in this plane.  

\begin{defin}
Given a pair of non-incident edges in a tetrahedron, the {\em medial 
parallelogram}
determined by the pair is the parallelogram whose vertices are the mid-points
of the remaining four edges.
\end{defin}

	Our main result is then

\begin{thm}
The area of the medial parallelogram determined by the edges of lengths
$d$ and $e$ in the tetrahedron of Figure \ref{tet} is

\[ \frac{1}{8}\sqrt{4d^2e^2 - (b^2 + f^2 - a^2 - c^2)^2} \]  

\end{thm}

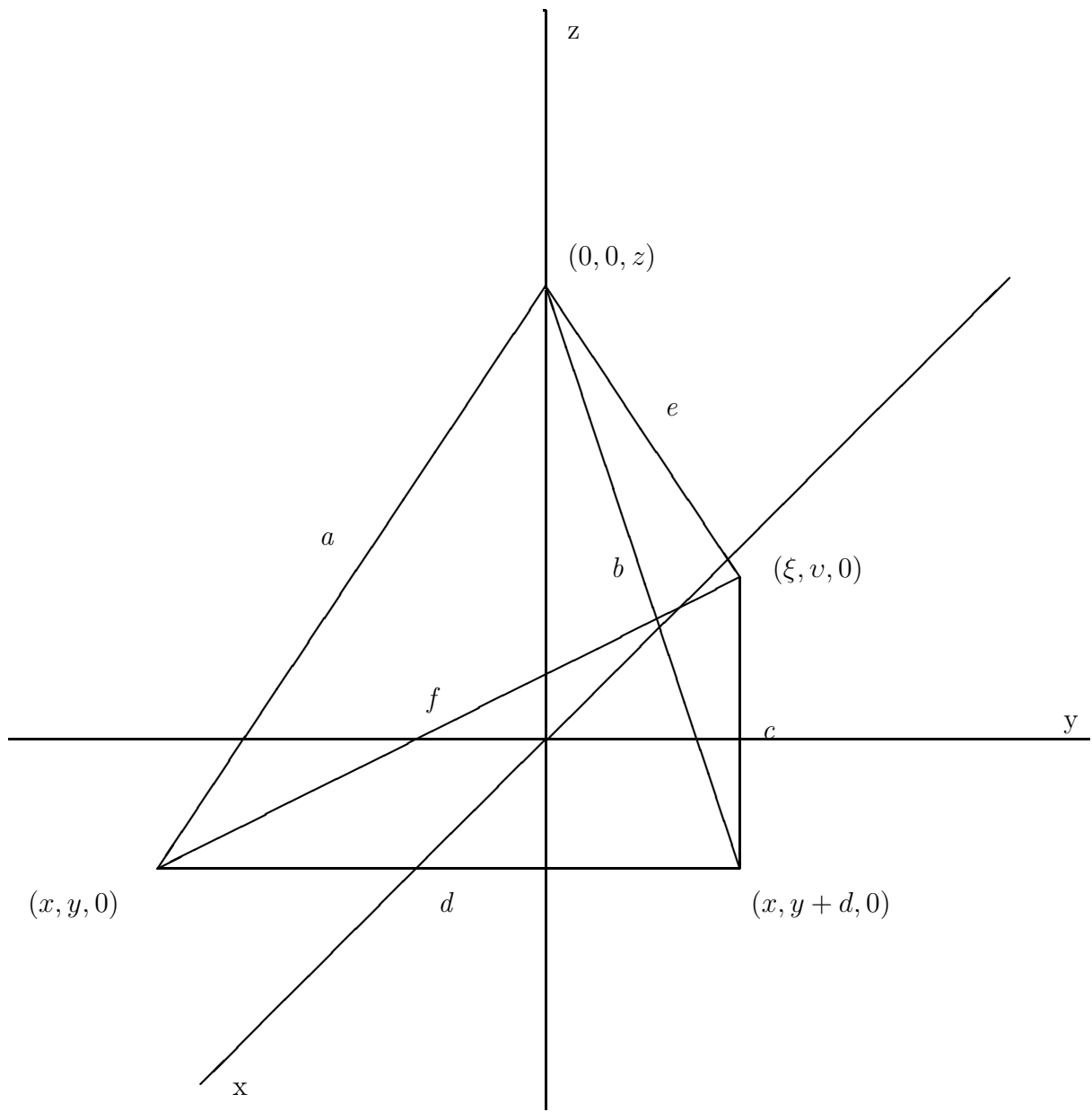
\begin{figure}

\setlength{\unitlength}{0.00083300in}%
\begingroup\makeatletter\ifx\SetFigFont\undefined
% extract first six characters in \fmtname
\def\x#1#2#3#4#5#6#7\relax{\def\x{#1#2#3#4#5#6}}%
\expandafter\x\fmtname xxxxxx\relax \def\y{splain}%
\ifx\x\y   % LaTeX or SliTeX?
\gdef\SetFigFont#1#2#3{%
  \ifnum #1<17\tiny\else \ifnum #1<20\small\else
  \ifnum #1<24\normalsize\else \ifnum #1<29\large\else
  \ifnum #1<34\Large\else \ifnum #1<41\LARGE\else
     \huge\fi\fi\fi\fi\fi\fi
  \csname #3\endcsname}%
\else
\gdef\SetFigFont#1#2#3{\begingroup
  \count@#1\relax \ifnum 25<\count@\count@25\fi
  \def\x{\endgroup\@setsize\SetFigFont{#2pt}}%
  \expandafter\x
    \csname \romannumeral\the\count@ pt\expandafter\endcsname
    \csname @\romannumeral\the\count@ pt\endcsname
  \csname #3\endcsname}%
\fi
\fi\endgroup
\begin{picture}(7539,7674)(2104,-7453)
\thicklines
\put(3151,-5761){\line( 0, 1){  0}}
\put(3151,-5761){\line( 1, 0){4050}}
\put(7201,-5761){\line(-1, 3){1350}}
\put(5851,-1711){\line(-2,-3){2700}}
\put(5851,-1711){\line( 2,-3){1350}}
\put(7201,-3736){\line(-2,-1){4050}}
\put(7201,-3736){\line( 0,-1){2025}}
\put(5851,-1741){\line( 0,-1){5700}}
\put(5851,-7441){\line( 0, 1){  0}}
\put(5851,-1726){\line( 0, 1){1935}}
\put(5851,209){\line(-1, 0){ 15}}
\put(2116,-4861){\line( 1, 0){7515}}
\put(5851,-4876){\line( 1, 1){3225}}
\put(5851,-4861){\line(-1,-1){2400}}
\put(4276,-3511){\makebox(0,0)[lb]{\smash{\SetFigFont{12}{14.4}{it}a}}}
\put(6301,-3736){\makebox(0,0)[lb]{\smash{\SetFigFont{12}{14.4}{it}b}}}
\put(7351,-4861){\makebox(0,0)[lb]{\smash{\SetFigFont{12}{14.4}{it}c}}}
\put(5101,-6061){\makebox(0,0)[lb]{\smash{\SetFigFont{12}{14.4}{it}d}}}
\put(6676,-2611){\makebox(0,0)[lb]{\smash{\SetFigFont{12}{14.4}{it}e}}}
\put(5026,-4636){\makebox(0,0)[lb]{\smash{\SetFigFont{12}{14.4}{it}f}}}
\put(6001,-1561){\makebox(0,0)[lb]{\smash{\SetFigFont{12}{14.4}{it}$(0,0,z)$}}}
\put(2251,-6061){\makebox(0,0)[lb]{\smash{\SetFigFont{12}{14.4}{it}$(x,y,0)$}}}
\put(7276,-6061){\makebox(0,0)[lb]
	{\smash{\SetFigFont{12}{14.4}{it}$(x,y+d,0)$}}}
\put(7426,-3736){\makebox(0,0)[lb]
	{\smash{\SetFigFont{12}{14.4}{it}$(\xi,\upsilon,0)$}}}
\put(3676,-7336){\makebox(0,0)[lb]{\smash{\SetFigFont{12}{14.4}{rm}x}}}
\put(9451,-4786){\makebox(0,0)[lb]{\smash{\SetFigFont{12}{14.4}{rm}y}}}
\put(6001, 14){\makebox(0,0)[lb]{\smash{\SetFigFont{12}{14.4}{rm}z}}}
\end{picture}

\caption{\label{coord.tet} The tetrahedron coordinatized}
\end{figure}

\noindent {\bf Proof:} The key is to set the tetrahedron in a three
dimensional coordinate system in a particularly convenient way.  We do
this by placing the vertex incident to the edges of lengths $a$, $b$, and
$e$ on the positive $z$-axis, and the remaining vertices in the $xy$-plane,
with the edge of length $d$ parallel to the $y$-axis.  The vertices
then have coordinates of the form $(0,0,z)$, $(x,y,0)$, $(x,y+d,0)$, and
$(\xi,\upsilon,0)$. (See Figure \ref{coord.tet}.)

	The medial tetrahedon is then spanned by the vectors 
$\vec{u} = (0,\frac{d}{2},0)$
and $\vec{v} = (-\frac{\xi}{2},-\frac{\upsilon}{2},\frac{z}{2})$, and thus

\begin{eqnarray}
{\rm Area} & = & \|\vec{u}\times \vec{v}\| \\
	& = & \left\| \left|
			\begin{array}{ccc}
			\vec{i} & \vec{j} & \vec{k} \\
			0 & \frac{d}{2} & 0 \\
			-\frac{\xi}{2} & -\frac{\upsilon}{2} & \frac{z}{2} 
			\end{array} \right| \right\| \\ 
	& = & \left\| \left( \frac{zd}{4},0,\frac{\xi d}{4} \right) \right\| \\
	& = & \frac{d}{4}\sqrt{z^2 + \xi^2}
\end{eqnarray}

	Now consider the relationships between the coordinates and the 
edge lengths:

\begin{eqnarray}
f^2 & = & (x-\xi)^2 + (y-\upsilon)^2 \\
c^2 & = & (x-\xi)^2 + (y+d-\upsilon)^2 \\
a^2 & = & x^2 + y^2 + z^2 \\
b^2 & = & x^2 + (y+d)^2 + z^2 \\
e^2 & = & \xi^2 + \upsilon^2 + z^2 
\end{eqnarray}

	Subtracting Equation (7) from Equation (8), and solving the resulting
equation for $y$ gives

\begin{eqnarray}
y & = & \frac{1}{2d}(b^2-a^2-d^2) .
\end{eqnarray}

	Subtracting Equation (5) from Equation (6), substituting the value
of $y$ from Equation (10), and solving for $\upsilon$ gives

\begin{eqnarray}
\upsilon & = & \frac{1}{2d}(b^2 + f^2 - a^2 - c^2)
\end{eqnarray}

	Now observe that by Equation (9) we have

\begin{eqnarray}
z^2 + \xi^2 & = & e^2 - \upsilon^2
\end{eqnarray}

	Substituting this, and the value for $\upsilon$ from Equation (11)
into Equation (4) and simplifying gives the desired result. $\Box$.

\end{document}